\documentclass{elsart3-1}

\usepackage{amssymb}

\usepackage[english,francais]{babel}

\newtheorem{theorem}{Theorem}[section]

\newtheorem{e-proposition}[theorem]{Proposition}

\newtheorem{e-definition}[theorem]{Definition\rm}

\newtheorem{theoreme}{Th\'eor\`eme}[section]
\newtheorem{lemme}[theoreme]{Lemme}
\newtheorem{proposition}[theoreme]{Proposition}
\newtheorem{corollaire}[theoreme]{Corollaire}

\renewcommand{\theequation}{\arabic{equation}}
\def\Fin#1{\leavevmode\unskip\nobreak\quad\hspace*{\fill}{#1}}

\newenvironment{preuve}{{\noindent} {\bf Preuve.}}{\Fin{$\square$}}

\newcommand{\diff}{{\rm d}}
\newcommand{\R}{{\mathbb R}}

\def\og{\leavevmode\raise.3ex\hbox{$\scriptscriptstyle\langle\!\langle$~}}
\def\fg{\leavevmode\raise.3ex\hbox{~$\!\scriptscriptstyle\,\rangle\!\rangle$}}

\begin{document}
\centerline{}
\begin{frontmatter}




%
\selectlanguage{francais}
\title{Forme semi-locale des feuilletages legendriens}



\author[saadi-zobida]{Sa\^adi Benabb\'es},
\ead{saadibenabbes$\underline{ }$fr@yahoo.fr}
\author[camille]{Camille Laurent-Gengoux}
\ead{camille.laurent-gengoux@univ-lorraine.fr}
\author[saadi-zobida]{Zobida Souici-Benhammadi}
\ead{zbenhamadi2000@yahoo.Fr}

\address[saadi-zobida]{D\'epart. de Math\'ematiques, Facult\'e des Sciences, Univsersit\'e Badji Mokhtar, Annaba, Alg\'erie}
\address[camille]{Institut Elie Cartan de Lorraine, Campus de Metz, France}



\begin{abstract}
\selectlanguage{francais}
Nous donnons une forme canonique semi-locale pour les feuilletages legendriens sur une vari\'et\'e de contact. Ce r\'esultat g\'en\'eralise une forme canonique locale donn\'ee par Libermann et Pang au voisinage d'une sous-vari\'et\'e de Legendre transverse, et est \`a mettre en parall\`ele avec un r\'esultat semi-local de Weinstein dans le cas symplectique. Au cours de la d\'emonstration, nous introduisons une classe de cohomologie qui mesure l'obstruction \`a rendre plat un feuilletage en modifiant la forme de contact.

\vspace{.5cm}

\selectlanguage{english}
\noindent{\bf Abstract}
\vskip 0.5\baselineskip
\noindent
{\bf Semi-local form of Legendrian foliations. }
We describe a semi-local canonical form for Legendrian foliations on contact manifolds in the neighbourhood of a Legendrian submanifold. This result generalizes local results by Libermann and Pang on Legendrian foliations on contact manifolds,  and is analogeous to a semi-local result by Weinstein in the symplectic case. For the proof, we introduce and use a class of cohomology that obstructs the possibility to make a Legendrian foliation flat.
\end{abstract}
\end{frontmatter}

\selectlanguage{english}
\section*{Abridged English version}

Weinstein \cite{W} proved that a neighbourhood of a Lagrangian submanifold transverse to a given Lagrangian foliation is symplectomorphic to the cotangent bundle of the transverse manifold. Under, this isomorphism, the leaves of the Lagrangian foliation correspond to the fiber the cotangent bundle. The cotagent bundle is therefore the unique semi-local model for Lagrangian foliations in a neighbourhood of a Lagrangian and transverse submanifold. For contact manifolds, what correspond to Lagrangian foliations (resp. submanifolds) are Legendrian foliations (resp. submanifolds). Libermann  \cite{Lib} and Pang  \cite{Pang} have established that, for a Legendrian submanifold transverse (in a sense given by (\ref{transverseEng}), which is not the usual sense) to a Legendrian foliation, the local model is the jet bundle. There is however, a function that comes into the picture and whose failure of being constant measures the non-flatness of the initial foliation. The purpose of the present note is to describe the semi-local model in this same very context.

Our work in mainly based on ideas and results by Libermann  \cite{Lib}, but it adds to them a new cohomological idea (which is in some sense implicit in \cite{Lib}). More precisely, we construct a class of cohomology that obstructs the possibility to gauge the contact form by multiplication by a non-zero function in order to make it flat - i.e. the Reeb vector has a flow that preserves the foliation. In the neighbourhood of a transversal Legendrian manifols, this obstruction class is zero, and therefore allows to gauge the contact form to find a simplified semi-local model, then to "un-gauge" it to find the semi-local model for the initial contact form. 

In short, there are two main results in this paper. The first one states that there is a class in the $H^1$
of functions constant on the leaves of a Legendrian foliation $ {\mathcal F}$ on a contact manifold $(M,\theta)$ which vanishes if and only if there is a non-zero function $F$ such that ${\mathcal F}$ is flat with respect to the contact form $F \theta$. The second one is the following:

\begin{theorem}
\label{theo:principal_English}
 Let $(M,\theta)$ be a contact manifold equipped with a Legendrian foliation $ {\mathcal F}$.
Any Legendrian submanifold $N$ of $M$ transverse to $ {\mathcal F}$ in the sense that, for all $n \in N $:
\begin{equation}\label{transverseEng} T_n N \oplus T_n {\mathcal F}\oplus \langle \eta \rangle = T_n M ,
\end{equation}
admits a neighbourhood diffeomorphic to a neighbourhood  of the zero section in the space of jets, namely
$T^*N \times {\mathbb R} \to  N$, through a diffeomorphisms that maps  $N$ to the zero section, that maps $ {\mathcal F}$ on the fibers of the projection $T^* N \times {\mathbb R} \to N \times {\mathbb R} $ and maps $\theta$ to 
a $1$-form that reads $H {\diff t} + \lambda $ with $H$ a function that never vanishes in a neighbourhood of the zero section and $ \lambda $ the Liouville form on $T^* N$.
\end{theorem}

\selectlanguage{francais}
\section{Introduction}
\label{}

Soit $ (M, \theta)$ une vari\'et\'e de contact \cite{Blair,LM} de dimension $2n+1$. On appelle \emph{feuilletage legendrien} un feuilletage ${\mathcal F}$ dont les feuilles sont des sous-vari\'et\'es de
dimension $n$ telles que, en tout point $m \in M$, l'espace tangent $T_m {\mathcal F} $ \`a la feuille passant par $m$ 
est compris dans le noyau de $\theta $. 

De m\^eme que toute vari\'et\'e symplectique de dimension $2n$ admet localement des coordonn\'ees de Darboux,
dont le mod\`ele local est un voisinage de la section nulle dans $T^* \R^n \to \R^n$, munie de la forme canonique
 $\sum_{i=1}^n \diff q_i \wedge \diff p_i$, il est bien connu \cite{Blair} que toute vari\'et\'e de contact a pour mod\`ele local le fibr\'e des jets, c'est-\`a-dire 
un voisinage de la section nulle dans $T^* {\R}^n \times \R \to \R^n \times \R$, 
munie de la forme canonique $\diff t + \sum_{i=1}^n p_i \diff q_i$. Dans le cas symplectique, lorsqu'un feuilletage lagrangien est donn\'e, ce diff\'eomorphisme local peut \^etre choisi de sorte que le feuilletage se transporte sur les fibres de la projection canonique $T^* \R^n \to \R^n$. Y. Pang \cite{Pang} et particuli\`erement P. Libermann \cite{Lib} d\'emontr\`erent que, dans le cas d'une vari\'et\'e de contact munie d'un feuilletage legendrien, 
ce feuilletage legendrien peut aussi \^etre transport\'e sur les fibres de la projection canonique
$T^* {\R}^n \times \R \to \R^n \times \R$. Dans ce cas n\'eanmoins, l'expression de la forme de contact ne peut pas \^etre totalement simplifi\'ee, mais on obtient seulement qu'elle est de la forme~:
 \begin{equation}\label{eq:formecanoniquelocale} H \diff t + \sum_{i=1}^n p_i \diff q_i \end{equation}
o\`u $H$ est une fonction qui d\'epend en g\'en\'eral de toutes les variables $ p_1,q_1, \dots,p_n,q_n, t$
et ne s'annule jamais dans un voisinage de la section nulle.

Dans le cas symplectique, ces th\'eor\`emes locaux ont \'et\'e g\'en\'eralis\'es en des th\'eor\`emes semi-locaux, c'est-\`a-dire  valables au voisinage d'une sous-vari\'et\'e $N$ lagrangienne transverse \`a un feuilletage lagrangien, par Weinstein \cite{W}. Dans cette note, nous proposons un analogue semi-local en g\'eom\'etrie de contact et au voisinage d'une sous-vari\'et\'e legendrienne $N$ tranverse (en un sens donn\'e dans l'\'equation (\ref{transverseFr}), sens qui diff\`ere du sens usuel)  au feuilletage legendrien.

Notre id\'ee est en fait de se ramener au cas plat, c'est-\`a-dire au cas o\`u la fonction $H$ qui appara\^it dans
(\ref{eq:formecanoniquelocale}) peut \^etre choisie constante - ce qui revient \`a supposer que le feuilletage legendrien est pr\'eserv\'e par le champ de Reeb. Un feuilletage legendrien $ {\mathcal F}$ n'est \'evidemment pas toujours plat, mais si on autorise \`a multiplier la forme de contact par une fonction qui ne s'annule jamais, on peut localement rendre ${\mathcal F} $ plat.  Une classe de cohomologie, qui se trouve dans une cohomologie naturellement associ\'ee au feuilletage, va donner une obstruction \`a la possibilit\'e de faire cette op\'eration globalement. Il se trouve qu'au voisinage d'une sous-vari\'et\'e de Legendre transverse, cette classe d'obstruction s'annule, ce qui permet la simplification.

Nous aimerions pr\'eciser aussi que notre r\'esultat  semi-global diff\`ere de celui obtenu par Pang \cite{Pang}
sur les feuilletages legendriens. Celui-ci en effet obtient des r\'esultats semi-globaux, voire globaux, mais en 
consid\'erant le cas o\`u les feuilles sont compactes et simplement connexes, ce qui ne peut jamais \^etre lorsqu'on se place au voisinage d'une transversale.


\section{Une obstruction \`a la platitude des feuilletages de Legendre}
\label{sec:applatissable}

Soit $(M ,\theta )$ une vari\'{e}t\'{e} de contact de dimension $2n+1$,
et $\eta$ son champ de Reeb. On munit $(M,\theta)$ d'un feuilletage de Legendre ${\mathcal F} $.
On dit que ${\mathcal F} $ est \emph{plat pour $\theta $} si
la forme fondamentale d\'efinie par Y. Pang \cite{Pang} ainsi~:
 $$  -\imath_\eta {\mathcal L}_X {\mathcal L}_Y \theta = 
 \imath_{[X,\eta] \wedge Y} \diff \theta  \quad \quad \hbox{ pour tous champs de vecteurs $X,Y$ tangents \`a ${\mathcal F}$}$$
est nulle (notons que cette notion correspond \`a la notion de feuilletage \emph{$\omega$-complet} chez \cite{Lib}).
Plus simplement, ainsi qu'il d\'ecoule\footnote{R\'esultat qui fut ensuite red\'emontr\'e dans \cite{Jayne}, lemme 4.2} des propositions 2.4, 2.6 et du corollaire 2.10 de P. Libermann \cite{Lib}, le feuilletage legendrien ${\mathcal F} $ est plat pour $\theta $ si et seulement si le champ de Reeb $\eta$ pr\'eserve $ {\mathcal F}$.

On dira qu'un feuilletage legendrien ${\mathcal F}$ sur une vari\'et\'e de contact $ (M,\theta)$ sur une vari\'et\'e de contact est \emph{applatissable} s'il existe une fonction $f \in {\mathcal C}^{\infty}(M,{\mathbb R}^*) $ telle que ${\mathcal F} $ soit plat pour la forme de contact $f \theta $. Notons que si une telle une fonction existe, on peut la choisir strictement positive. Nous allons dans cette section caract\'eriser en terme cohomologiques les feuilletages de Legendre applatissables.

\begin{lemme} \label{lem:troispoints}
Soit $N$ une vari\'et\'e et soit $H $ une fonction, d\'efinie dans un voisinage de $N $ dans  $T^*N \times {\mathbb R} $, qui ne s'annule jamais. Il  existe un voisinage de $N  $ dans  $T^*N \times {\mathbb R} $ sur lequel~:
\begin{enumerate}
 \item $H \diff t + \lambda$ est une forme de contact (ici $ \lambda $ d\'esigne la forme de Liouville sur $T^* N$),
 \item l'application $m_H$ de $ T^*N \times \R \to T^*N \times \R $ d\'efinie par
 $$ (\alpha,t) \mapsto ( H\alpha, t)  \hbox{ pour tous $ \alpha \in T^*M$, $ t \in {\mathbb R}$,} $$
 est un diff\'eomorphisme sur son image,
 \item Ce diff\'eomorphisme \'echange les formes de contact $ H (\diff t + \lambda ) $ et  $((m_H^{-1})^* H ) \diff t + \lambda$.
\end{enumerate}
\end{lemme}
\begin{preuve}
Soit $x_1,\dots,x_n)$ un syst\`eme de coordonn\'ees locales sur ${\mathcal U} \subset M$, $p_1, \dots,p_n$ les coordonn\'ees duales sur le $ T^*{\mathcal U}$ de sorte que  $(x_1, \dots,x_n,p_1, \dots,p_n,t) $ est un syst\`eme de coordonn\'ees sur ${\mathcal U} \times \R $. Le premier point de ce lemme se d\'emontre par un calcul direct en coordonn\'ees locales: on peut aussi utiliser le th\'eor\`eme 3.3 dans \cite{Lib}. Le second point vient de ce que, dans ces coordonn\'ees, $m_H$ s'\'ecrit~:
 $$ (x_1, \dots,x_n,p_1, \dots,p_n,t) \mapsto (x_1, \dots,x_n,H(x,p,t) \, p_1, \dots,H(x,p,t) \, p_n,t) ,$$
 et de ce que la diff\'erentielle de cette application est un multiple de l'identit\'e en tout point o\`u $p_1=\dots=p_n=0 $, le coefficient de multiplication \'etant $ H(x,0,t)$. Le dernier point vient de ce que le tir\'e en arri\`ere de 
 $H(x,p,t) (\diff t + \sum_{i=1}^n p_i \diff x_i) $ par $ m_H$ est $ H( x,\frac{p}{H(x,p,t)}, t ) \diff t +  \sum_{i=1}^n p_i \diff x_i$, comme un calcul direct le montre.
\end{preuve}

\begin{prop} \label{prop:feuilletage_plat}
Soit $(M,\theta)$ une vari\'et\'e de contact.
\begin{enumerate}
 \item Tout feuilletage de Legendre est localement applatissable.
 \item Soit ${\mathcal F} $ un feuilletage de Legendre plat pour $\theta $ et  $f\in {\mathcal C}^{\infty}(M,{\mathbb R}^*) $. Le feuilletage $ {\mathcal F}$ est plat pour $f \theta$ si et seulement si $f $ est constant sur les feuilles du feuilletage $ {\mathcal F}$.
\end{enumerate}
\end{prop}
\begin{preuve}
 Le premier point de cette proposition d\'ecoule du troisi\`eme point du lemme \ref{lem:troispoints}
 et du th\'eor\`eme 3.3 dans \cite{Lib}. Pour le second point, 
 il est connu (voir par exemple Marle-Libermann \cite{LM}) que pour toute fonction $f \in C^*(M,{\mathbb R}_+^*)$, $f \theta$ est encore une forme de contact et que son champ de Reeb est $ \frac{\eta }{f}+{\mathcal X}_{{\rm ln}(f)}$
o\`{u} ${\mathcal X}_{{\rm ln}(f)}$ est le champ hamiltonien de ${\rm ln}(f)$.
\end{preuve}

Nous allons construire une classe dans le premier espace de cohomologie du faisceau des fonctions constantes sur les feuilles du feuilletage ${\mathcal F} $, faisceau que l'on notera encore ${\mathcal F} $, dont les sections sur un ouvert ${\mathcal U} \subset M$ sont not\'ees ${\mathcal F}_{\mathcal U} $.

Par le premier point du lemme \ref{prop:feuilletage_plat}, il existe un recouvrement $({\mathcal U}_i)_{i\in I} $ de $M$ et, pour chaque indice $i \in I$, une fonction $f_i \in {\mathcal C}^{\infty}({\mathcal U}_i,{\mathbb R}_+)$ telle que $ \left.{\mathcal F}\right|_{{\mathcal U}_i}$ soit plat pour $f_i \theta $.
Quitte \`a remplacer $f_i$ par $-f_i$, on peut supposer que ces fonctions sont strictement positives.

Pour toute paire d'indice $(i,j) \in I^2$, la fonction $ f_{ij} := f_i/f_j$ est telle que ${\mathcal F} $ est \`a la fois plat pour la forme de contact
$f_i \theta $ et pour la forme de contact $ f_j \theta = f_{ji} f_i \theta $. Le second point du lemme \ref{prop:feuilletage_plat} implique que 
 $ f_{ij} \in {\mathcal F}_{{\mathcal U}_{ij}}$, et donc que $ {\rm ln}(f_{ij}) \in  {\mathcal F}_{{\mathcal U}_{ij}}$. 
 Cette derni\`ere famille de fonctions v\'erifie \'evidemment~:
  $$ {\rm ln}(f_{ij})+ {\rm ln}(f_{jk})+ {\rm ln}( f_{ki})
  = {\rm ln}(f_{i})- {\rm ln}(f_{j})+ {\rm ln}(f_{j})- {\rm ln}(f_{k})+{\rm ln}(f_{k})- {\rm ln}(f_{i})
  =0$$
et donc d\'efinit une classe dans le premier espace de cohomologie $ H^1({\mathcal F})$ du faisceau $ {\mathcal F}$, classe que l'on notera $ [Obs(\theta,{\mathcal F})]$. 

\begin{theoreme}\label{theo:applatissable}
Soit $ (M,\theta,{\mathcal F})$ une vari\'et\'e de contact munie d'un feuilletage legendrien ${\mathcal F} $.
La classe $ [Obs(\theta,{\mathcal F})]\in H^1({\mathcal F})$ est bien d\'efinie\footnote{Autrement dit, ne d\'epend
pas des choix faits lors de sa construction} et v\'erifie  $ [Obs(\theta,{\mathcal F})]= [Obs(F\theta,{\mathcal F})] $
pour toute fonction $F \in {\mathcal C}^{\infty}(M,\R^*)$. De plus, les points suivants sont \'equivalents:
\begin{enumerate}
 \item[(i)] $ [Obs(\theta,{\mathcal F})]=0$,
 \item[(ii)] le feuilletage ${\mathcal F} $ est applatissable.
 \end{enumerate}
\end{theoreme}
\begin{preuve}
 Le second point du lemme \ref{prop:feuilletage_plat} implique qu'un autre choix des fonctions $f_i$ ci-dessus serait de la forme $h_i f_i $, o\`u les fonctions $h_i$ appartiennent \`a ${\mathcal F}_{{\mathcal U}_i} $, ce qui revient \`a remplacer ${\rm ln}(f_{ij})$
 par ${\rm ln}(f_{ij}) + {\rm ln}(h_{i}) - {\rm ln}(h_{j}) $, autrement dit \`a ajouter un cobord. 
 
 
 De plus $ [Obs(\theta,{\mathcal F})]=0$ si et seulement si il existe
 un recouvrement ${\mathcal U}_i $ et des fonctions $h_i$ appartenant \`a ${\mathcal F}_{{\mathcal U}_i} $
 telles que ${\rm ln}(f_{ij}) + {\rm ln}(h_{i}){\rm ln}(h_{j}) =0$. Pour tout indice $i$, la forme $h_i \theta $
 est une forme de contact globalement d\'efinie sur $ M$ et le feuilletage est plat pour cette forme.
 La r\'eciproque est \'evidente.
\end{preuve}

Dans le cas des feuilletages de Legendre ${\mathcal F} $, il y a une notion de \emph{sous-vari\'et\'e transverse}  qui consiste\footnote{La terminologie est ici trompeuse : une vari\'et\'e transverse \`a un feuilletage de Legendre n'est pas transverse au feuilletage.} \`a exiger que pour tout $n $ dans la vari\'et\'e transverse $ N$:
\begin{equation}\label{eq:transverse}
T_n N \oplus T_n {\mathcal F}_n \oplus \langle \eta \rangle = T_n M
\end{equation}
 Des consid\'erations cohomologiques g\'en\'erales et le th\'eor\`eme \ref{theo:applatissable}  impliquent le r\'esultat suivant.

\begin{corollaire}\label{coro:voisinagetranverse}
Soit $ (M,\theta,{\mathcal F})$ une vari\'et\'e de contact munie d'un feuilletage de Legendre ${\mathcal F} $.
Toute sous-vari\'et\'e $N$  transverse au feuilletage (au sens donn\'e par (\ref{eq:transverse})) 
admet un voisinage  ${\mathcal V} $ auquel la restriction du feuilletage ${\mathcal F} $ est applatissable.
\end{corollaire}
\begin{preuve}
Ce corollaire vient tout simplement de ce qu'il existe un voisinage ${\mathcal V}$ de $N$
sur lequel $H^1({\mathcal F})=0$. Ceci vient \`a son tour de ce qu'il existe un voisinage de $ N \times \{0\}$
dans $N \times \R$ tel que l'image de l'application $N \times \R \to M$ d\'efinie par $(n,t) = \varphi_t (n)$
($\varphi_t$ \'etant le flot du champ de Reeb) est une sous-vari\'et\'e $N'$ transverse (au sens usuel) au feuilletage 
$ {\mathcal F}$ et donc un voisinage ${\mathcal V} $ de $N'$ dans $M$ sur lequel le feuilletage $ {\mathcal F}$ d\'efinit une fibration sur $N'$ dont les fibres sont, topologiquement, des boules ouvertes.
Sur cet ouvert, on a \'evidemment un isomorphisme entre $H^1(M,{\mathcal F})$ et le $H^1$
du faisceau des fonctions lisses sur $N$, lequel est nul. Le th\'eor\`eme \ref{theo:applatissable} permet alors de conclure.
\end{preuve}

\section{Forme semi-locale des feuilletage de Legendre plats}

Nous donnons maintenant un r\'esultat semi-local dans le cas plat.

\begin{proposition}\label{prop:lengendreplat}
 Soit $(M,\theta)$ une vari\'et\'e de contact munie d'un feuilletage legendrien plat. Toute vari\'et\'e legendrienne $N$ transverse\footnote{Ce qui signifie, on le rappelle, que l'\'equation (\ref{eq:transverse}) est satisfaite, et diff\`ere du sens usuel.} au feuilletage ${\mathcal F} $ admet un voisinage diff\'eomorphe \`a un voisinage de la section nulle dans l'espace des jets $T^*N \times {\mathbb R}$ par un diff\'eomorphisme qui
\begin{enumerate}
\item \'echange $ N$ et la section nulle de $ T^*N \times {\mathbb R} \to N$,
\item \'echange $ {\mathcal F}$ et les fibres de la projection $T^* N \times {\mathbb R} \to N \times {\mathbb R} $ 
\item \'echange la forme $\theta$ et la $1$-forme $ {\diff t} + \lambda $,  o\`u $ \lambda $ est la forme de Liouville sur $T^* N$.
\end{enumerate}
 \end{proposition}
\begin{preuve}
Quitte \`a se restreindre \`a un voisinage ${\mathcal V} $ de $N$ dans $M$, on peut supposer que l'union $N_{\mathcal F}$ des feuilles  de $ \left. {\mathcal F}\right|_{\mathcal V}$ qui 
intersectent $N$ est une sous-vari\'et\'e (de dimension $2n$) de $ {\mathcal V} \subset M$.
La forme $\diff \theta$, restreinte \`a $N_{\mathcal F}$, est symplectique.
Pour cette forme  $ \omega:= \left. \diff \theta\right|_{N_{\mathcal F}}$, le feuilletage $ {\mathcal F}$ est un feuilletage lagrangien et $N$ est une sous-vari\'et\'e lagrangienne. Par un th\'eor\`eme de Weinstein, il existe un voisinage de $p(N)$ dans $ N_{\mathcal F}$ isomorphe, via un diff\'eomorphisme $\Psi$, \`a un voisinage de la section nulle dans $T^*N $, diff\'eomorphisme qui \'echange $\omega$ et la $2$-forme canonique du cotangent, que l'on note encore $\omega $, et qui
\'echange la restriction de $ {\mathcal F}$ \`a $N_{\mathcal F} $ avec les fibres de la projection $T^*N \to N $.

Comme le champ de Reeb est transverse \`a la sous-vari\'et\'e $ {\mathcal N}_{\mathcal F} $ et conserve le feuilletage,
l'application $ \Phi: T^* N \times \R \to {\mathcal V} $ d\'efinie par~:
 $$ (\alpha, t) \mapsto \phi_t (\Psi(\alpha) )  .$$
o\`u $\phi_t$ est le flot du champ de Reeb, est un diff\'eomorphisme d'un voisinage ${\mathcal W} $ de 
la section nulle dans l'espace des jets $T^*N \times {\mathbb R} \to N$ dans un voisinage de $N$ dans $M$ que l'on appellera encore ${\mathcal V} $.

Montrons que ce diff\'eomorphisme $\Phi$ convient ce pour quoi il suffit de d\'emontrer que $\theta'= \Phi^* \theta$ 
co\"incide avec $ \diff t + \lambda$, les conditions \emph{(i)} et \emph{(ii)} \'etant \'evidemment v\'erifi\'ees. 

Commen\c{c}ons par \'etablir quelques propri\'et\'es. Comme $ \Phi$ \'echange le champ de Reeb $ \eta$ avec le champ de vecteur $\frac{\partial }{\partial t} $, cedernier est lle champ de Reeb de $\theta'$, soit~:
\begin{equation}
 \label{eq:theta'd/dt}
\imath_{\frac{\partial }{\partial t} } \theta' =1 \mbox{ et } {\mathcal L}_{\frac{\partial }{\partial t}} \theta' = 0 =0
\end{equation}
Comme la restriction \`a $ T^*N \times \{0\} $ de $ \diff \theta'$ co\"incide avec la forme symplectique canonique de cet espace, les relations (\ref{eq:theta'd/dt}) impliquent la relation suivante~:
\begin{equation}
 \label{eq:omega}\diff \theta' = \diff \lambda .
 \end{equation}
Comme la section nulle et les fibres de la projection sont des sous-vari\'et\'es legendriennes, on a aussi~:
\begin{equation}
 \label{eq:theta'feuilles} \left. \theta' \right|_N=0 \mbox{ et } \left. \theta'\right|_F=0 \hbox{ pour toute fibre de la projection canonique.} 
 \end{equation}
   
	Toutes ces propri\'et'es impliquent que $\theta' $ co\"incide avec $ \diff t + \lambda$ et ach\`eve donc la preuve de \emph{(iii)}. Pour commencer, la partie gauche de (\ref{eq:theta'd/dt}) implique que $ \theta' = \diff t + \mu_t$ o\`u $\mu_t$ est, pour tout $t$, une $1$-forme sur $T^*N$. La partie droite de (\ref{eq:theta'd/dt})  implique que cette $1$-forme ne d\'epend en fait pas de $t$. On la note par $\mu$. Par (\ref{eq:omega}), la $1$-forme $\theta' - \lambda - \diff t = \mu - \lambda$ v\'erifie $ \diff (\lambda -\mu) =0$. Ceci implique que, localement,  $\diff (\mu -\lambda) = \diff g $ pour une certaine fonction $g$ qui ne d\'epend pas de $t$. La  partie gauche de (\ref{eq:theta'feuilles}) implique que
    $g$ est constante le long des fibres de la projection canonique, et provient donc d'une fonction $h$ sur $N$. La  partie droite de (\ref{eq:theta'feuilles}) implique alors que cette fonction $h$ est \`a son tour constante, et donc que $\mu = \lambda $, ce qui implique \emph{(iii)} et ach\`eve la d\'emonstration du th\'eor\`eme.
\end{preuve}

\section{Forme semi-locale des feuilletage de Legendre~: cas g\'en\'eral}

Soit $(M,\theta,{\mathcal F})$ une vari\'et\'e de contact munie d'un feuilletage legendrien, et $N$ une sous-vari\'et\'e de Legendre transverse au feuilletage. Le corollaire \ref{coro:voisinagetranverse} implique qu'il existe une fonction $f  \in {\mathcal C}^{\infty}(M,{\mathbb R}_+^*)$
telle que ${\mathcal F} $ est plat pour $f\theta $, pourvu que l'on se restreigne \`a un voisinage 
suffisamment petit de $N$. Par la proposition \ref{prop:lengendreplat}, il existe un diff\'eomorphisme 
$\Phi$ de vari\'et\'es de Lengendre entre un voisinage ${\mathcal U}$ de $N$ dans $M$ et un voisinage de la section nulle 
$(T^*N \times \R) \to N $ qui envoie $ {\mathcal F}$ sur le feuilletage naturel $ \tilde{\mathcal F}$ 
et $ f\theta$ sur la forme de contact naturelle $\diff t + \lambda $.
En cons\'equence, $ \Phi$ envoie la forme de contact $\theta $ sur $F (\diff t + \lambda ) $, o\`u $ F =\Phi^* f$ est une fonction qui ne s'annule jamais sur un voisinage de la setion nulle. Le lemme \ref{lem:troispoints} donne alors le th\'eor\`eme suivant~:

\begin{theoreme}
\label{theo:principal}
 Soit $(M,\theta)$ une vari\'et\'e de contact munie d'un feuilletage legendrien. Toute vari\'et\'e legendrienne $N$ transverse\footnote{Ce qui signifie, on le rappelle, que l'\'equation (\ref{eq:transverse}) est satisfaite, et diff\`ere du sens usuel.} au feuilletage ${\mathcal F} $ admet un voisinage diff\'eomorphe \`a un voisinage de la section nulle dans l'espace des jets $T^*N \times {\mathbb R} \to N$ par un diff\'eomorphisme qui 
\begin{enumerate}
\item \'echange $ N$ et la section nulle de $ T^*N \times {\mathbb R} \to N$,
\item \'echange $ {\mathcal F}$ et les fibres de la projection $T^* N \times {\mathbb R} \to N \times {\mathbb R} $ 
\item \'echange la forme $\theta$ et la $1$-forme $ H {\diff t} + \lambda $, o\`u $H$ est une fonction sur $T^*N \times {\mathbb R}$ qui ne s'annule jamais sur un voisinage de
la section nulle et o\`u $ \lambda $ est la forme de Liouville sur $T^* N$.
\end{enumerate}
\end{theoreme}




\end{document}